# AN ANALYTICAL GOMBOC
M. L. Sloan

**Abstract:**

Investigation of the mathematical requirements for a three dimensional geometrical object to qualify as a Gomboc (mono-monostatic) has resulted in the discovery of two specific, analytical Gomboc shapes. Analytical in that the function describing the Gomboc surface is infinitely differentiable.

In this brief note, the analysis undertaken is summarized and the formulae for the two specific shapes provided.

## 1. Introduction

A Gomboc is a three dimensional convex solid of constant density which possesses only two equilibrium points, one stable and the other, unstable, when resting on a horizontal surface in the presence of a constant vertical gravitational field.

The question of the existence of such an object was raised by Russian mathematician Vladimir Arnold in 1995, with the problem solved in 2006 by Gábor Domokos and Péter Várkonyi of Hungary (P.L. Várkonyi and G. Domokos: Mono-monostatic bodies: the answer to Arnold's question. The Mathematical Intelligencer, Volume 28, Number 4., pp 34-38. (2006)) [1].

The solution found by Domokos and Várkonyi required extreme precision in any physical embodiment, with tolerances of $10^{-5}$ in their earliest embodiment, and somewhat more reasonable, but still extreme, tolerances of $10^{-3}$ in later embodiments. [As discussed in various internet articles on the Gomboc *e.g. en.wikipedia.org* and *plus.maths.org* ]

In this brief note, two specific analytic Gomboc shapes are presented which are simple in form and, to the extent examined, not as sensitive to variations.

This presentation will be expository, with none of the details of the multiple, often fruitless, though instructive, research paths undertaken.

## 1. Overall Gomboc Requirements

In general, the bounding surface of an analytic Gomboc may specified by a function:

$\psi(z, x, y) = $ **constant**

where **x**, **y**, **z** are standard Cartesian coordinates describing the 3 dimensional space in which the Gomboc is imbedded.

We orient the coordinate system so that the origin (**x = y = z = 0**) is the center of mass of the homogeneous density Gomboc.

In examining the requirements for a Gomboc, polar coordinates **r**, $\theta$, $\phi$ are appropriate, with:

$x = r \sin(\theta) \cos(\phi)$
$y = r \sin(\theta) \sin(\phi)$
$z = r \cos(\theta)$

The center of mass of the gomboc is then the point **r = 0**.

In polar coordinates, the bounding surface $\psi$ **(z, x, y) = constant** may then be inverted to yield a solution in **r** for the boundary of the Gomboc:

$r = F(\theta, \phi)$

The requirement of convexity ensures that **r** is single valued in $\theta$, $\phi$ and moreover positive definite.

### Center of Mass Requirements

Since **r = 0** is the center of mass of the Gomboc, we require:

$\int d^3V \; x(r, \theta, \phi) = 0$

$\int d^3V \; y(r, \theta, \phi) = 0$       where $d^3V = r^2 \sin(\theta) \, dr \, d\theta \, d\phi$, with the integration carried out over the volume of the Gomboc         (1)

$\int d^3V \; z(r, \theta, \phi) = 0$

The boundary **r = F ( θ , φ)** being single valued in θ , φ , the **r** integration is trivial and one arrives at the follow center of mass requirements:

$$\int_0^{2\pi} d\phi \int_0^{\pi} \sin(\theta) \cos(\theta) \, F^4(\theta, \phi) \, d\theta = 0$$

$$\int_0^{2\pi} d\phi \int_0^{\pi} \sin^2(\theta) \sin(\phi) \, F^4(\theta, \phi) \, d\theta = 0 \qquad (2)$$

$$\int_0^{2\pi} d\phi \int_0^{\pi} \sin^2(\theta) \cos(\phi) \, F^4(\theta, \phi) \, d\theta = 0$$

Using complex variable notation: **exp( i φ ) = cos(φ) + i sin(φ)** , the second two equations may be combined to allow a more succinct statement of the center of mass requirements:

$$\int_0^{2\pi} d\phi \int_0^{\pi} \sin(\theta) \cos(\theta) \, F^4(\theta, \phi) \, d\theta = 0 \qquad (3a)$$

$$\int_0^{2\pi} d\phi \int_0^{\pi} \sin^2(\theta) \exp(i\phi) \, F^4(\theta, \phi) \, d\theta = 0 \qquad (3b)$$

### Equilibria Requirements

Given a Gomboc shape defined by:

$$r - F(\theta, \phi) = 0 \qquad (4)$$

equilibrium points are those points on the surface where the normal to the surface

$$\nabla ( r - F(\theta, \phi) ) \qquad (5)$$

lies wholly in the $\hat{r}$ unit normal direction.

Specifically, then, equilibria are those θ , φ sets of points where the components of the gradient ∇ in the θ and φ direction vanish:

$$1/ \sin(\theta) \; \partial_\phi F(\theta, \phi) = 0 \qquad (6)$$

and

$$\partial_\theta F(\theta, \phi) = 0 . \qquad (7)$$

For a shape to be a Gomboc, there can only be two sets of such points.

### Convexity Requirement

As with all convex Gomboc shapes, the radius of curvature must be positive at each point on the Gomboc surface.  For the specific analytic Gombocs exhibited below, that is easily achieved.

## 2.  Specific Analytic Gombocs

Minimally "bumpy" shapes probably stand the best chance of meeting the Gomboc requirements, since Gombocs are required to exhibit two and only two equilibrium points.  That observation, along with many hours of research and false starts, have led the author to restrict investigations of analytic Gomboc geometries to the simplest $\phi$ dependency possible [1].

$$F^4(\theta, \phi) = R(\theta) + \sin(\theta) A(\theta) \cos(\phi - P(\theta)) , \; A(\theta) > 0 \qquad (8)$$

Among such possible solutions, the following particular embodiment has proven particularly useful:

$$r^4 = F^4(\theta, \phi) = 1 + 4\beta \sin(\theta) \cos(\phi - P(\theta)) , \qquad (9)$$

with $\beta$ a small positive constant.[2]

---

[1]   It is well know that a purely $\phi$-symmetric Gomboc is not possible.

[2]   This Gomboc shape is essentially a sphere of unit radius with small surface perturbations. To lowest order in $\beta$, this particular embodiment may be functionally realized by slicing a unit sphere transverse to the **z** axis into many (infinitesimally) thin disks and then re-assembling the sphere from the disks, but with the center of each disk offset from the **z** axis by $(x_o, y_o) = (\beta \cos P(\theta), \beta \sin P(\theta))$.

This formulation greatly simplifies the equilibria and center of mass requirements.

In particular, there are two and only two equilibrium points:

$$\theta = \pi/2, \quad \phi = P(\pi/2) + \pi \quad \text{(Stable equilibrium point)}$$

$$\text{and} \quad \theta = \pi/2, \quad \phi = P(\pi/2) \quad \text{(Unstable equilibrium point)}.$$

(10)

Moreover, the center of mass requirement reduces to the restriction on $P(\theta)$:

$$\int_0^\pi \sin^3(\theta) \exp(i P(\theta)) \, d\theta = 0 \quad (11)$$

Thus, Eq. (9) properly describes an analytic Gomboc, provided $P(\theta)$ satisfies the center of mass requirement Eq. (11).

Two particular examples of Analytic Gombocs are now presented:

### Analytic Gomboc 1

The choice of $P(\theta) = n\,\theta$ satisfies Eq. (11) for all $n = 2p+1$, with $p \geq 2$.

Choosing $p = 2$ ($n = 5$) to provide the simplest and least "bumpy" solution, we arrive at one embodiment of an analytic Gomboc:

$$r^4 = 1 + 4\beta \sin(\theta) \cos(\phi - 5\theta), \quad (12)$$

with a single stable equilibrium point at $\theta = \pi/2, \phi = 3\pi/2$

(13)

and a single unstable equilibrium point at $\theta = \pi/2, \phi = \pi/2$.

The positive constant $\beta$ must, of course, be sufficiently small to satisfy the convexity requirement. $\beta = 0.15$ or smaller appears to suffice.

### Analytic Gomboc 2

Introducing the variable $\eta = 3\pi/2 \, (\cos(\theta) - \cos^3(\theta)/3)$
and expressing $P(\theta)$ in terms of $\eta$, the center of mass requirement, Eq. (11), reduces to the simple requirement:

$$\int_{-\pi}^{\pi} \exp(i\, P(\eta))\, d\eta = 0 \qquad (14)$$

The following particularly simple solution:

$$p(\eta) = \eta = 3\pi/2 \,(\cos(\theta) - \cos^3(\theta)/3) \qquad (15)$$

provides a second analytic Gomboc solution:

$$r^4 = 1 + 4\beta \sin(\theta) \cos(\phi - 3\pi/2\,(\cos(\theta) - \cos^3(\theta)/3)) \qquad (16)$$

with a single stable equilibrium point at   $\theta = \pi/2, \phi = \pi$

(17)

and a single unstable equilibrium point at  $\theta = \pi/2, \phi = 0$.

Examination of the convexity requirement for this Gomboc shape indicates that a value of  $\beta = 0.17$  or smaller should suffice.

Note that the Eq. (12) Gomboc wraps smoothly two and one-half revolutions around the **z** axis ( $\phi = 0$ to $5\pi$ as $\theta$ traverses **0** to **π** ), while the second Eq (16) Gomboc exhibits only one full revolution around the **z** axis, but with a somewhat nonlinear $\phi$ vs $\theta$ path.

Finally, as indicated in Footnote 2, the solutions presented amount to surface perturbations on a unit sphere.  Accordingly, the solutions may be scaled up by any constant  $r_o^4$  to achieve any specific size Gomboc desired.

**ACKNOWLEDGEMENT:**  Posting and dissemination of this technical note would not have been possible without the guidance of Professor Gabor Domokos, co-inventor of the Gomboc, whose support and encouragement are greatly appreciated.